\documentclass[11pt]{article}

\usepackage{amsmath}
\usepackage{amssymb}

\newcommand{\C}{{\mathbb  C}}
\newcommand{\D}{{\mathbb D}}
\newcommand{\R}{{\mathbb  R}}

\renewcommand{\P}{{\mathbb  P}}

\newcommand{\T}{{\mathbb  T}}
\newcommand{\A}{{\cal A}}
\newcommand{\F}{{\cal F}}

\newcommand{\OO}{{\cal O}}
\newcommand{\PSH}{{\operatorname{{\cal PSH}}}}
\newcommand{\SH}{{\operatorname{{\cal SH}}}}
\newcommand{\sing}{{\operatorname{sing}}}
\renewcommand{\phi}{\varphi}

\newtheorem{theorem+}           {Theorem}      [section]
\newtheorem{definition+}  [theorem+]  {Definition}
\newtheorem{lemma+}  [theorem+]  {Lemma}
\newtheorem{corollary+}  [theorem+]  {Corollary}
\newtheorem{proposition+}  [theorem+]  {Proposition}
\newtheorem{example+}  [theorem+]  {Example}
\newtheorem{question+}  [theorem+]  {Question}

\newenvironment{theorem}{\begin{theorem+}\sl}{\end{theorem+}\rm}
\newenvironment{definition}{\begin{definition+}\rm}{\end{definition+}\rm}
\newenvironment{lemma}{\begin{lemma+}\sl}{\end{lemma+}\rm}

\newenvironment{proposition}{\begin{proposition+}\sl}{\end{proposition+}\rm}

\newenvironment{proof}{\medbreak\noindent{\it Proof:}\rm}{\hfill$\square$\rm}
\newenvironment{prooftx}[1]{\medbreak\noindent{\it #1:}\rm}{\hfill$\square$\rm}

\begin{document}

\title{%
EXTREMAL $\omega$-PLURISUBHARMONIC FUNCTIONS AS ENVELOPES OF DISC FUNCTIONALS
}

\author{%
Benedikt Steinar Magn\'usson
}

\maketitle

\begin{abstract}
For each closed, positive $(1,1)$-current $\omega$ on a complex manifold $X$ and
each $\omega$-upper semicontinuous function $\phi$ on $X$ we associate a
disc functional and prove that its envelope is equal to the supremum of
all $\omega$-plurisubharmonic functions dominated by $\phi$.
This is done by reducing to the case where 
$\omega$ has a global potential. 
Then the result follows from Poletsky's theorem,
which is the special case $\omega=0$. Applications of this result include a
formula for the relative extremal function of an open set in $X$ and, in some cases,
a description of the $\omega$-polynomial hull of a set.
\end{abstract}


\section{Introduction}\label{intro}
If $\omega$ is a closed, positive $(1,1)$-current on a connected complex manifold $X$, 
then for every point $x_0 \in X$ we can find a neighbourhood
$U$ of $x_0$ and a plurisubharmonic local potential $\psi$ for $\omega$, 
i.e., $dd^c \psi = \omega$ on $U$. 
Let $u\colon X  \to \overline \R$ be a function on $X$ with values in the extended real line.
If we can locally write $u = v - \psi$, where $v$ is plurisubharmonic, then we say the 
function $u$ is \emph{$\omega$-plurisubharmonic}.
We denote by $\PSH(X,\omega)$ the set of all $\omega$-plurisubharmonic 
functions on $X$ which are not identically equal to $-\infty$ in any connected
component of $X$.

If $\psi_1$ and $\psi_2$ are two local potentials for $\omega$ then their difference is 
pluriharmonic on their common set of definition. This implies that the singular 
set, \emph{$\sing(\omega)$}, of $\omega$ is well defined and locally given as 
$\psi^{-1}(\{-\infty\})$ for a local potential $\psi$ of $\omega$.

We say that a function $\phi\colon X \to \overline \R$ is \emph{$\omega$-upper semicontinuous} 
if $\phi + \psi$ is upper semicontinuous on $U \setminus \sing(\omega)$, extends 
to an upper semicontinuous function on $U$ for every local potential 
$\psi \colon U \to \R \cup \{-\infty\}$ of $\omega$, and for $a\in \sing(\omega)$
we have $\limsup_{X\setminus \sing(\omega) \ni z \to a} u(z)=u(a)$.

An \emph{analytic disc} is a holomorphic map $f\colon \D \to X$ from the unit disc $\D$ into $X$. 
It is said to be \emph{closed} if it can be 
extended to a holomorphic map in some neighbourhood of the closed unit disc. 
We let $\mathcal O(\D,X)$
denote the set of all analytic discs 
and $\A_X$ denote the set of all closed analytic discs in $X$.

For every analytic disc we have a pullback $f^*\omega$ of $\omega$ which is a Borel-measure on $\D$. 
It is defined locally as the Laplacian of the pullback $f^*\psi$, of a local potential $\psi$, 
to an open subset of $\D$. We define $R_{f^*\omega}$
as the Riesz potential of this measure on $\D$. 

The main result of this paper is the following
\begin{theorem}\label{th} Let $X$ be a connected complex manifold, 
$\omega$ be a closed, positive $(1,1)$-current on $X$, and $\phi$ be an 
$\omega$-upper semicontinuous function on $X$ such that 
$\{u \in \PSH(X,\omega) ; u \leq \phi \}$ is nonempty. 
Then for $x \in X\setminus \sing(\omega)$
$$
\sup\{u(x) ; u \in \PSH(X,\omega), u \leq \phi \} = 
\inf\{ -R_{f^*\omega}(0) + \int_\T \phi \circ f\, d\sigma ; f \in \A_X, f(0) = x \},
$$
where $\sigma$ is the arc length measure on the unit circle $\T$ normalized to $1$. 
Furthermore, if $\{ u \in \PSH(X,\omega)  ; u \leq \phi \}$ is empty then the right
hand side is $-\infty$.
\end{theorem}

This theorem is a generalization of Poletsky's theorem, which is the special case $\omega=0$, 
see Poletsky \cite{Pol:1991}, 
L\'{a}russon and Sigurdsson \cite{LarSig:1998} and \cite{LarSig:2003}, and 
Rosay \cite{Ros:2003}.

However, if $\omega$ has a global potential $\psi$, i.e.~$\psi \in \PSH(X)$ 
with $dd^c \psi = \omega$, then the formula 
above becomes 
$$
\sup\{u(x) ; u \in \PSH(X,\omega), u \leq \phi \} + \psi(x) = 
\inf\{ \int_\T (\psi + \phi) \circ f\, d\sigma  ; f \in \A_X, f(0) = x \},
$$
which is a direct consequence of Poletsky's theorem since $\psi + \phi$ is 
an upper semicontinuous function.
This case is handled in Theorem \ref{th_stein}.

The general case follows from this case and an $\omega$-version of a
reduction theorem (Theorem 1.2 in \cite{LarSig:2003}) proved by L\'arusson and Sigurdsson, 
see Theorem \ref{th_red}.
The reduction theorem states that 
the right hand side in Theorem 1.1 is $\omega$-plurisubharmonic on $X$ 
if all its pullbacks to a manifold with a global potential are, and if we can assume some 
continuity properties of it
with respect to the discs in $\A_X$.

By applying Theorem \ref{th} to the characteristic function of the complement of an open set $E$ we get a disc formula
for the relative extremal function, which Guedj and Zeriahi introduce in \cite{GueZer:2005}, Chapter 4. 
Our result is the following 
\begin{multline*}
\sup\{ u(x)  ; u \in \PSH(X,\omega), u|_E \leq 0 \text{ and } u\leq 1 \}\\
= \inf\{ -R_{f^*\omega}(0) + \sigma(\T \setminus f^{-1}(E) ) ; f \in \A_X, f(0) = x \}.
\end{multline*}
In certain cases this formula can give us a description of the $\omega$-polynomial hull of a set, which 
is a generalization of the polynomial hull in $\C^n$.

For more information about the recent development of $\omega$-psh
functions we refer the reader to 
Guedj and Zeriahi \cite{GueZer:2005} and \cite{GueZer:2007}, 
Harvey and Lawson \cite{HarLaw:2006}, 
Ko{\l}odziej \cite{Kol:2003},
Dinew \cite{Din:2007}, and
Branker and Stawiska \cite{BraSta:2009}. 
In these papers $X$ is usually assumed to be a compact {K\"ahler} manifold and $\omega$ a 
smooth current on $X$.

{\bf Acknowledgement:}
I would like to thank my teacher Professor Ragnar Sigurdsson for pointing out
this subject to me and for helping me writing the paper. 
I would also like to thank Professor J\'on Ing\'olfur Magn\'usson for his help and Institut Mittag-Leffler 
where a part of this
work was done during the Several Complex Variables program in 2008.


\section{Basic results on $\omega$-plurisubharmonicity}\label{sec:2}
Here we will define $\omega$-plurisubharmonic functions and study their properties analogous to those
of plurisubharmonic functions.

Assume $X$ is a complex manifold of dimension $n$ and is $\omega$ a closed positive $(1,1)$-current 
on $X$, i.e. $\omega$ acts on $(n-1,n-1)$-forms. 

It follows from Proposition  1.19, Ch.~III in \cite{Demailly}, 
that locally there is a plurisubharmonic function $\psi$ such that $dd^c \psi = \omega$. 
Here  $d$ and $d^c$ are the
real differential operators $d = \partial + \overline \partial$ and $d^c = i(\overline \partial - \partial)$.
Hence, in $\C$ then $dd^c u = \Delta u\, dV$, where $\Delta u$ is the Laplacian of $u$
and $dV$ is the standard volume form.

Note that the difference of two potentials for $\omega$ is a pluriharmonic function, thus $\mathcal C^\infty$. 
This implies that the singular set of $\omega$, $\sing(\omega)$, is well defined as
the union of all $\psi^{-1}(\{-\infty\})$, for all local potentials $\psi$ of $\omega$.

In the case when $\omega$ has continuous local potentials we have no trouble with continuity. If $\psi$ is a continuous local
potential for $\omega$ then $u + \psi$ is upper semicontinuous if and only if $u$ is. In general
this is not always the case, and we do not want to exclude the case when $\psi$ takes the value 
$-\infty$, e.g., when $\omega$ is a current of integration, or if $\psi$
is discontinuous.
This however forces us to define the value of $u+\psi$ at points $x \in \sing(\omega)$ where
$\psi(x)=-\infty$ and possibly $u(x)=+\infty$. If $u+\psi$ is bounded above on 
$X\setminus \sing(\omega)$ in a neighbourhood of $x$, then this can be done by taking 
upper limits of $u+\psi$ as we approach points in $\sing(\omega)$,
we therefore make the following definition.

\begin{definition} 
A function $u\colon X \to [-\infty,+\infty]$ is called \emph{$\omega$-upper semicontinuous} ($\omega$-usc) 
if for every 
$a \in \sing(\omega)$, $\limsup_{X\setminus \sing(\omega) \ni z\to a} u(z) = u(a)$ and for 
each local potential 
$\psi$ of $\omega$, defined on an open subset $U$ of $X$,
$u+\psi$ is upper semicontinuous on $U\setminus \sing(\omega)$, and locally bounded above 
around each point of $\sing(\omega)$.
\end{definition}

Equivalently we could say that $\limsup_{X\setminus \sing(\omega) \ni z\to a} u(z) = u(a)$ for every $a \in \sing(\omega)$ and
$u+\psi$ extends as
$$
 \limsup_{U\setminus \sing(\omega) \ni z\to a} (u+\psi)(z), \qquad \text{for } a \in \sing(\omega)
$$ 
to an upper semicontinuous function on $U$ with values in
$\R \cup \{-\infty\}$. We will denote this extension by $(u+\psi)^\star$.

Note that the question whether $(u+\psi)^\star$ is usc does not depend at all on the values of 
$u$ at $\sing(\omega)$.
The reason for the conditions on $u$ at $\sing(\omega)$ is to ensure that
$u$ is Borel measurable and to uniquely determine the function from its
values outside of $\sing(\omega)$.

It is easy to see that $u$ is Borel measurable from the fact that
$u = (u+\psi)-\psi$ is the difference of two Borel measurable functions 
on $X\setminus \sing(\omega)$ and that $u$ restricted to the Borel set $\sing(\omega)$ is 
the increasing limit of usc functions. Hence it is Borel measurable.

\begin{definition}
A function $u\colon X\to [-\infty,+\infty]$ is called 
\emph{$\omega$-pluri\-subharmonic} ($\omega$-psh) if it is $\omega$-usc and
$u+\psi$ is psh on $U\setminus \sing(\omega)$ for every local potential $\psi$ of $\omega$
defined on an open subset $U$ of $X$.
We let $\PSH(X,\omega)$ denote the set of all $\omega$-psh functions on $X$ 
which are not identically equal to $-\infty$ on any connected component of $X$.
When the manifold is one dimensional we say that these functions are
$\omega$-subharmonic ($\omega$-sh) and denote the set of $\PSH(X,\omega)$ by $\mathcal{SH}(X,\omega)$.
\end{definition}

Equivalently we could say that $\limsup_{X\setminus \sing(\omega) \ni z\to a} u(z) = u(a)$ 
for every $a \in \sing(\omega)$ and
$(u+\psi)^\star$ is psh for every local potential $\psi$. 
We see that the $\omega$-psh functions are locally integrable because outside of the zero
set $\sing(\omega)$ they can locally be written as the difference of two functions which 
are locally integrable on $X$.

Our approach depends on the fact that we can define the pullback of currents by holomorphic maps. 
This we can do in two very different cases,
first if the map is a submersion and secondly if it is an analytic disc not lying in 
$\sing(\omega)$. 

If $\Phi \colon Y \to X$ is a submersion and $\omega$ is a current on $X$ then we can define the inverse
image $\Phi^*\omega$ of $\omega$ by its action on forms, 
$\langle \Phi^*\omega,\tau \rangle = \langle \omega, \Phi_*\tau \rangle$, where 
$\Phi_*\tau$ is the direct image of the form $\tau$. For more details
see Demailly \cite{Demailly}, 2.C.2 Ch.~I.

If $f$ is an analytic disc in $X$, $f(\D) \nsubseteq \sing(\omega)$ 
then we can define a closed, positive $(1,1)$ current $f^*\omega$ on $\D$
in the following way.

Let $a \in \D$ and $\psi \colon U \to \R \cup \{-\infty\}$ be a local potential on an open neighbourhood 
$U$ of $f(a)$, and let $V$ be the connected component of $f^{-1}(U)$ containing $a$. 
If $\psi \circ f \neq -\infty$ on $V$ then we define
$f^*\omega = dd^c (\psi \circ f)$. 
If $\psi \circ f = -\infty$ on $V$ then we define $f^*\omega$ as the  
measure which sends
$\emptyset$ to $0$ and $E$ to $+\infty$ for all $E \neq \emptyset$. We denote this measure by $+\infty$.
Observe that if $\psi \circ f = -\infty$ on $V$ then the same applies for every other local potential.  
In fact, if $\psi$ and $\psi'$ are two potentials for $\omega$ on open sets $U$ and $U'$ respectively 
then on the open set $V \cap V'$ we have two subharmonic functions $f^*\psi$ and $f^*\psi'$
which differ by a harmonic function. Therefore, if one of them is equal to $-\infty$ the other one 
is also equal to $-\infty$.

Remember that every positive $(n,n)$-current (of order 0) on an $n$ 
dimensional manifold can be given by 
a positive Radon measure, and conversely every positive Radon measure defines a 
positive $(n,n)$-current.
So when we pull $\omega$ back to $\D$ by an analytic disc it is possible look at 
it both as a $(1,1)$-current and as a Radon measure.

We let $R_{f^*\omega}$ be the Riesz potential of the positive measure $f^*\omega$. It is defined by 
\begin{equation}\label{R_f}
R_{f^*\omega}(z) = \int_\D G_\D(z,\cdot)\, d(f^*\omega),
\end{equation}
where $G_\D$ is the Green function for the unit disc, $G_\D(z,w) = \log (|z-w|/|1-z\overline w|)$.
If $f^*\omega = +\infty$ then $R_{f^*\omega} = -\infty$. The Riesz potential of $f^*\omega$ is not 
identically $-\infty$
if and only if $f^*\omega$ satisfies the Blaschke condition (see \cite{Hor:convexity}, Theorem 3.3.6)
$$
\int_\D (1-|\zeta|)\, d(f^*\omega)(\zeta) < +\infty.
$$

If $f$ is a closed analytic disc not lying in $\sing(\omega)$, then this condition is satisfied since $f^*\omega$ is a Radon measure in a neighbourhood of the unit disc, thus with finite mass on $\D$.

Also if we have a local potential $\psi$ defined in a neighbourhood of $\overline{f(\D)}$ then the Riesz
representation formula, (ibid.), at the point 0 gives
\begin{equation}\label{riesz}
 \psi(f(0)) = R_{f^*\omega}(0) + \int_\T \psi \circ f\, d\sigma.
\end{equation}

\begin{proposition}\label{equi}
 The following are equivalent for a 
 function $u$ on $X$.
\begin{enumerate}
 \item[(i)]\ $u$ is in $\PSH(X,\omega)$.
 \item[(ii)]\ $u$ is $\omega$-usc and $h^*u \in \SH(\D,h^*\omega)$ for all $h\in \A_X$ such that 
	 $h(\D) \nsubseteq \sing(\omega)$.
 \end{enumerate}
\end{proposition}
\begin{proof} 
Assume $u \in \PSH(X,\omega)$, take $h \in \A_X$, $h(\D) \nsubseteq \sing(\omega)$, and 
$a \in \D$. Let $\psi$ be a local potential 
for $\omega$ defined in a neighbourhood $U$ of $h(a)$. 
Note that $(u+\psi)^\star \circ h = (u\circ h + \psi \circ h)^\star$, i.e.~the extension of 
$(u + \psi) \circ h$ over $\sing(h^*\omega)$ is the same as the
extension of $u+\psi$ over $\sing(\omega)$ pulled back by $h$, for both functions are subharmonic and
equal almost everywhere, thus the same.
Since $(u + \psi)^\star \in \PSH(U)$ and
$(u + \psi)^\star \circ h = (u \circ h + \psi \circ h)^\star$ is subharmonic in a neighbourhood of $a$
we see that $u \circ h$ is $h^*\omega$-subharmonic.

Assume now that \textsl{(ii)} holds and let $\psi \in \PSH(U)$ be a local potential for $\omega$.
Then $(u+\psi)^\star$ is upper semicontinuous, and \textsl{(ii)} implies that $(u+\psi)^\star \circ h \in \SH(\D)$ for
\end{proof}

From the proposition we also see that $\omega$-plurisubharmonicity like plurisubharmonicity is a 
local property, so in condition \textsl{(ii)} it is sufficient to look at $h \in \A_U$ in a neighbourhood $U$ of a given point.

If $u_0 \neq -\infty$ is psh and $\phi$ is an usc function then the family 
$\{ u \in \PSH(X)  ; u_0 \leq u \leq \phi \}$ is compact in the $L^1_\text{loc}$ topology, 
which implies that $\sup\{ u(x)  ; u \in \PSH(X), u \leq \phi \}$ is plurisubharmonic.
We have a similar result for $\omega$-plurisubharmonic functions.

\begin{proposition}
If $\phi\colon X \to [-\infty,+\infty]$ is $\omega$-usc, 
 $\mathcal F_{\omega,\phi} = \{ u \in \PSH(X,\omega) ; u \leq \phi \}$ and
 $\mathcal F_{\omega,\phi} \neq \emptyset$, then 
 $\sup \mathcal F_{\omega,\phi} \in \PSH(X,\omega)$,
and consequently $\sup \mathcal F_{\omega,\phi} \in \mathcal F_{\omega,\phi}$.
\end{proposition}
\begin{proof}
Define $s_\phi = \sup\mathcal F_{\omega,\phi}$, by definition $s_\phi + \psi \leq \phi + \psi$ 
outside of $\sing(\omega)$ for
every local potential $\psi$ on $U$. Since $\phi$ is $\omega$-usc,
the upper semicontinuous regularization $(s_\phi+\psi)^*$of $s_\phi +\psi$ also satisfies
$$
(s_\phi + \psi)^* \leq \phi + \psi \text{ on } U\setminus \sing(\omega).
$$
Note that the left hand side is plurisubharmonic. 
We define the function $S$ on $X \setminus \sing(\omega)$ by
$$
S(x) = (s_\phi + \psi)^*(x) - \psi(x),
$$
where $\psi$ is a local potential for $\omega$ in some neighbourhood of $x$. Observe that
since the difference between two local potentials is continuous it is clear that the function
$S$ is well defined on $X \setminus \sing(\omega)$. 
We extend $S$ to a $\omega$-usc function on $X$ by taking the $\limsup$ at points in $\sing(\omega)$.
Furthermore, it is then obvious that $S \in \PSH(X,\omega)$ and 
$s_\phi \leq S \leq \phi$ follows from the inequality above, so $s_\phi = S \in \mathcal F_{\omega,\phi}$.
\end{proof}


\section{Disc functionals and their envelopes}\label{sec:3}
A \emph{disc functional} $H$ is a function defined on some subset $\A$ of $\OO(\D,X)$, the set
of all analytic discs in a manifold $X$, with values in $[-\infty,+\infty]$.
The \emph{envelope} $EH$ of a disc functional $H$ is then
a function defined on the set $X_\A = \{ x\in X  ; x=f(0) \text{ for some } f \in \A \}$
by the formula
$$
	EH(x) = \inf\{ H(f) ;  f\in \A, f(0)=x \}, \quad x \in X_\A.
$$

The Poisson disc functional is defined as $H_\phi(f) = \int_\T \phi\circ f\, d\sigma$ where $\phi$
is an usc function and $\sigma$ is the normalized arc length measure on the unit circle $\T$.
If $f \in \A_X$ is a closed analytic disc, $u \in \PSH(X)$ and $u \leq \phi$ then
$$
u(f(0)) \leq \int_\T u \circ f\, d\sigma \leq H_\phi(f).
$$
The envelope $EH_\phi$ of $H_\phi$ is plurisubharmonic by the Poletsky theorem \cite{LarSig:2003} and 
equal to the supremum of all plurisubharmonic functions less than or equal $\phi$,
that is
$$
\sup \{ u(x) ; u \in \PSH(X), u \leq \phi \} = \inf \Big\{ \int_\T \phi \circ f\, d\sigma ; f \in \A_X, f(0) = x \Big\}.
$$

We will now generalize the definition of the Poisson functional to $\omega$-usc functions and 
look at the largest $\omega$-psh minorant of $\phi$. This functional will be denoted
by $H_{\omega,\phi}$.

Fix a $\omega$-usc function $\phi$ on a complex manifold $X$ 
and a point $x\in X \setminus \sing(\omega)$, 
let $f \in \A_X$, $f(0)=x$ and assume there is a 
function $u \in \PSH(X,\omega)$, $u \leq \phi$.
Then since $f(0)=x\notin \sing(\omega)$ the pullback 
$f^*\omega$ is a well defined Radon measure on $\D$.
Remember that $R_{f^*\omega}$ is a global potential for 
$f^*\omega$ on $\D$ and equal to $0$ on the
boundary, so by Proposition \ref{equi},
$$
u(x) + R_{f^*\omega}(0) \leq \int_\T u\circ f\, d\sigma + \int_\T R_{f^*\omega}\, d\sigma.
$$
Since $u \leq \phi$ and $R_{f^*\omega} = 0$ on $\T$ we see that
\begin{equation}\label{fund_inequal}
 u(x) \leq - R_{f^*\omega}(0) + \int_\T \phi \circ f\, d\sigma.
\end{equation}
The right hand side is independent of $u$, so we define the functional $H_{\omega,\phi}$
for every $f\in \A_X$, $f(0)\notin \sing(\omega)$ by
$$
H_{\omega,\phi}(f) =  - R_{f^*\omega}(0) + \int_\T \phi \circ f\, d\sigma.
$$

Now take supremum on the left hand side of (\ref{fund_inequal}) over all $\omega$-psh function $u$ satisfying 
$u\leq \phi$ and infimum
on the right over all $f \in \A_X$ such that $f(0)=x$. Then we get the fundamental inequality
$$
\sup \F_{\omega,\phi} \leq EH_{\omega,\phi} \qquad \text{on } X\setminus \sing(\omega).
$$
Theorem \ref{th} states that this is actually an equality.

Since the disc functional is not defined for discs centered at $a \in\sing(\omega)$ we extend 
the envelope to a function on the whole space $X$ by  
$$
EH_{\omega,\phi}(a) = \limsup_{X\setminus \sing(\omega) \ni z \to a} EH_{\omega,\phi}(z).
$$

In the following we let $D_r = \{ t \in \C ; |t|<r\}$ and if $x$ is a point in $X$ then 
$H(x)$ will denote the value of $H$ at the constant disc $t \mapsto x$, the meaning should
always be clear from the context. 

Notice now that if we look at the constant discs in $X$, then $H_{\omega,\phi}(x) = \phi(x)$ 
and consequently $EH_{\omega,\phi} \leq \phi$.
Therefore, if we show that $EH_{\omega,\phi}$ is $\omega$-psh then it is in $\F_{\omega,\phi}$ and we have
an equality $\sup \F_{\omega,\phi} = EH_{\omega,\phi}$.	

\medskip

An immediate corollary of the main theorem is a formula for the \emph{relative extremal function} of a set $E$ in $\Omega$, 
where $\Omega$ is an open subset of $X$. It is defined as
$$
h_{E,\Omega,\omega}(x) = \sup\{u(x) ; u \in \PSH(\Omega,\omega), u|_E \leq 0 \text{ and } u\leq 1 \},
$$
Now assume $E$ is open and apply Theorem \ref{th} to $\Omega$ with $\phi$ as the characteristic function for the complement of $E$.
For $x\in \Omega \setminus \sing(\omega)$ it gives that
$$
h_{E,\Omega,\omega}(x) = \inf\{ -R_{f^*\omega}(0) + \sigma(\T \setminus f^{-1}(E)) ; f \in \A_\Omega, f(0) = x \}.
$$
When $\Omega = X$ we denote this function by $h_{E,\omega}$.

\medskip

In the local theory, $\omega=0$ and $X \subset \C^n$, the relative extremal function can be
used to describe the polynomial hull of a compact set.
The result is due to Poletsky \cite{Pol:1993}
and can also be found in the following form in \cite{LarSig:2007}, Theorem 2.
It states that for a compact set $K$ in $\C^n$, 
a point $a$ in $\C^n$, and $\Omega$ a pseudoconvex neighbourhood of $K$ and $a$, bounded
and Runge, the following are equivalent.
\begin{enumerate}
	\item[\textsl{(i)}]\ $a$ is in the polynomial hull of $K$.
	\item[\textsl{(ii)}]\ For every neighbourhood $U$
		of $K$ and every $\epsilon > 0$ there is a $f \in \A_\Omega$ with $f(0)=a$ and
		$\sigma(\T \setminus f^{-1}(U)) < \epsilon$.
\end{enumerate}
If we now wish to use our formula to get a similar result on a general complex manifold we have to start by
finding an alternative to the polynomial hull.
It follows from Theorem 5.1.7 in \cite{Kli} that the polynomial hull in 
$\C^n$ is exactly
the hull with respect to the psh functions in $\C^n$ of logarithmic 
growth.
These functions correspond to the $\omega$-psh functions on $\P^n$ if
$\omega$ is the integration current for the hyperplane at infinity.
This motivates the following definition, which is similar to the definition given 
by Guedj \cite{Gue:1999} of the $\omega$-polynomial hull.

\begin{definition}If $K \subset X$ is compact subset of a complex 
manifold $X$ and $\omega$ is a closed, positive
$(1,1)$-current on $X$, we define the $\omega$-polynomial hull of $K$ as
$$
\hat K^\omega = \{ x \in X ; u(x) \leq \sup_K u \text{ for all } u \in \PSH(X,\omega)\}.
$$
\end{definition}

Our goal is to use the disc formula above to describe this hull. To make that possible
we have to be able to use the relative extremal function to describe the hull, that
is there is an $\Omega \subset X$ such that $h^{-1}_{K,\Omega,\omega}(\{0\}) = \hat K^\omega$. 
In the local theory it is sufficient to have $\Omega \subset \C^n$ hyperconvex.

The disc formula only applies to open sets, here we are considering 
compact sets so we start by showing that it is enough to look at shrinking neighbourhoods
of compact sets.

\begin{proposition}\label{compact_limit}
Assume $\omega$ is a closed, positive $(1,1)$-current on a complex manifold $X$ such that
$\omega$ has continuous local potentials on $X\setminus \sing(\omega)$.
Let $K_1 \supset K_2 \supset ...$ be sequence of compact subsets of an open set
$\Omega \subset X \setminus \sing(\omega)$ and $K = \cap_{j=1}^\infty K_j$, then
$$
\lim_{j\to  \infty} h_{K_j,\Omega,\omega} = h_{K,\Omega,\omega}.
$$
\end{proposition}
The proof is the same as in the case $\omega=0$, see Klimek \cite{Kli} Proposition 4.5.10.
We only have to note that the assumptions on $\omega$
imply that all $\omega$-psh functions are usc.

\medskip

Next we derive the result with some assumptions on $\Omega$, below we see
that in some cases we can take $\Omega=X$.

\begin{proposition}
Let $K$ be a compact subset of $\Omega \subset X\setminus \sing(\omega)$, and assume $\omega$ has
continuous local potentials and $\Omega$ satisfies
$h^{-1}_{K,\Omega,\omega}(\{0\}) = \hat K^\omega$.
Then a point $x \in \Omega$ is in $\hat K^\omega$ if and only if for every neighbourhood $U$ of $K$ in
$\Omega$ and every $\epsilon >0$ there is an analytic disc $f \in \A_\Omega$ such that
$f(0)=x$ and 

$$
-R_{f^*\omega}(0) + \sigma(\T \setminus f^{-1}(U) ) < \epsilon.
$$
\end{proposition}
\begin{proof}
	Let $x \in \hat K^\omega$. Then $0 \leq h_{U,\Omega,\omega}(x) \leq h_{K,\Omega,\omega}(x)=0$, and by the
	disc formula for $h_{U,\Omega,\omega}$ there is a disc $f\in \A_\Omega$ such that $f(0)=x$ and
	$$
		-R_{f^*\omega}(0) + \sigma(\T \setminus f^{-1}(U) ) < \epsilon.
	$$

	Conversely, if such $f$'s exist then $h_{U,\Omega,\omega}(x) = 0$ for every neighbourhood $U$ of 
	$K$. Let $\{K_j\}$ be a sequence of compact subsets of $\Omega$ such that $\cap_j K_j = K$ and
	$K_{j+1} \subset \mathring K_j$. Then $h_{K_j,\Omega,\omega}(x)=0$ and by Proposition 
	\ref{compact_limit} we see that $h_{K,\Omega,\omega}(x) = 0$.
\end{proof}

\begin{proposition}
	If $\sing(\omega) =\emptyset$ and $X$ is compact then $h^{-1}_{K,\omega}(\{0\}) = \hat K^\omega$.
\end{proposition}
\begin{proof}
	Assume $x \in h^{-1}_{K,\omega}(\{0\})$ and let $u \in \PSH(X,\omega)$. 
	Note that if we let $\psi_j\colon U_j \to \R^+$ be positive local potentials for $\omega$, such that $\cup_j U_j = X$,
	then for every $j$ the function $u+\psi_j$ is usc and locally bounded on $U_j$, and by a compactness
	argument we then see that $\sup_X u < \infty$. 
	
	Now let $\tilde u= (u-\sup_K u)/(\sup_X u - \sup_K u)$ if $\sup_X u -\sup_K u > 1$, 
	else if $\sup_X u - \sup_K u \leq 1$
	then we let $\tilde u = u - \sup_K u$. Either way, $\tilde u$ is $\omega$-psh, $\tilde u \leq 1$ and $\tilde u|_K \leq 0$.
	Therefore $\tilde u \leq h_{K,\omega}$ and $\tilde u(x) \leq h_{K,\omega}(x) = 0$, that is $u(x) \leq \sup_K u$ and $x \in \hat K^\omega$.

	Obviously $\hat K^\omega \subset h^{-1}_{K,\omega}(\{0\})$, so $h^{-1}_{K,\omega}(\{0\}) = \hat K^\omega$.
\end{proof}


\section{Proof of the main result}\label{sec:4}

We start by restricting to the case when $\omega$ has a global potential.
The general case then follows from the reduction theorem later on.

\begin{lemma}\label{global_potential}
Let $\omega$ be a closed positive $(1,1)$-current on Stein manifold $X$. If there
is a current $\eta$ such that $d\eta = \omega$, then $\omega$ has a global
plurisubharmonic potential $\psi \colon X \to \R \cup \{-\infty\}$, 
so in particular $dd^c \psi = \omega$.
\end{lemma}
\begin{proof}
Since $\omega$ is a positive current it is real, 
so $\eta$ can be assumed to be real, $\eta \in \Lambda'_1(X,\R)$. 
Now write 
$\eta = \eta^{1,0} + \eta^{0,1}$, where $\eta^{1,0} \in \Lambda'_{1,0}(X,\C)$
and $\eta^{0,1} \in \Lambda'_{0,1}(X,\C)$. 
Note that $\eta^{0,1} = \overline{ \eta^{1,0}}$ since $\eta$ is real.
We see, by counting degrees, that 
$\overline \partial \eta^{0,1} = \omega^{0,2} = 0$, then since $X$ is Stein 
there is a distribution $\mu$ on $X$ such that 
$\overline \partial \mu = \eta^{0,1}$.
Then
$$
\eta = \overline{ \overline\partial \mu} + \overline\partial \mu
= \partial \overline \mu + \overline \partial \mu.
$$
If we set $\psi = (\mu - \overline \mu)/2i$, then
$$
\omega = d\eta = d(\partial \overline \mu + \overline \partial \mu)
= (\partial + \overline \partial)(\partial \overline \mu 
+ \overline \partial \mu)
= \partial \overline \partial (\mu - \overline \mu) 
= dd^c \psi.
$$
Finally, $\psi$ is a plurisubharmonic
function since $\omega$ is positive.
\end{proof}

\begin{theorem}\label{th_stein}
Let $\omega$ be a closed, positive $(1,1)$-current on a manifold $X$ and 
$\phi\colon X \to [-\infty,+\infty]$ an $\omega$-usc function such that 
$\mathcal F_{\omega,\phi} \neq \emptyset$.
If $\omega$ has a global potential $\psi$ then $EH_{\omega,\phi} \in \PSH(X,\omega)$ 
and consequently $EH_{\omega,\phi} = \sup \mathcal F_{\omega,\phi}$ 
on $X\setminus \sing(\omega)$.
\end{theorem}

\begin{proof}
For $f \in \A_X$, $f(\D) \nsubseteq \sing(\omega)$, the Riesz representation (\ref{riesz}) of $f^*\psi$ gives the following
$$
H_{\omega,\phi}(f) + \psi(f(0)) = 
H_{\omega,\phi}(f) + R_{f^*\omega}(0) + \int_\T \psi \circ f \, d\sigma = 
\int_\T (\psi + \phi)^\star \circ f\, d\sigma = H_{\psi+\phi}(f).
$$
The equality in the middle follows from the fact that 
$\int_\T \phi \circ f\, d\sigma + \int_\T \psi \circ f\, d\sigma = \int_\T (\phi+\psi)^\star\circ f\, d\sigma$
since
$\sigma(f^{-1}(\sing({\omega}))\cap \T) = 0$.
Therefore 
$$
EH_{\omega,\phi}(x) + \psi(x) = \inf\{ H_{\omega,\phi}(f) + \psi(x) ; 
f\in \A_X, f(0) = x \} = EH_{(\psi+\phi)^\star}(x).
$$
By Poletsky's theorem $EH_{(\psi+\phi)^\star}$ is psh, hence 
$EH_{\omega,\phi}$ is $\omega$-psh.
\end{proof}

If $\Phi\colon Y \to X$ is a holomorphic map between complex manifolds and $H$ is a disc functional 
on $X$, then the pullback $\Phi^*H$ is a disc functional on $Y$ defined by
$\Phi^*H(f) = H(\Phi\circ f)$ for $f\in \A_Y$. 
Since $\{ \Phi \circ f ; f\in \A_Y \} \subset \A_X$ we get the following.

\begin{lemma}\label{pullback}
 $\Phi^*EH \leq E\Phi^*H$ and equality holds if every disc in $Y$ is a lifting of a disc in $X$ by $\Phi$.
\end{lemma}

Moreover, for the Poisson functional $H_{\omega,\phi}$ we have the following result.

\begin{lemma}\label{submersion}
Assume $\Phi\colon Y \to X$ is a holomorphic submersion, then $\Phi^*H_{\omega,\phi} = H_{\Phi^*\omega,\Phi^*\phi}$.
\end{lemma}
\begin{proof}
By associativity of compositions we have $(\Phi_*f)^*\omega = f^*(\Phi^*\omega)$ for $f\in \A_Y$, 
$f(\D) \nsubseteq \Phi^{-1}(\sing(\omega))$, so
 \begin{eqnarray*}
\Phi^*H_{\omega,\phi}(f) = H_{\omega,\phi}(\Phi_* f) &=&
 -R_{(\Phi_* f)^*\omega}(0) + \int_\T \phi \circ \Phi \circ f\, d\sigma \\
&=& -R_{f^*(\Phi^*\omega)}(0) + \int_\T (\Phi^* \phi) \circ f\, d\sigma = H_{\Phi^*\omega,\Phi^*\phi}(f)
\end{eqnarray*}
\end{proof}

We will now state the reduction theorem which will enable us to prove Theorem \ref{th} using 
Theorem \ref{th_stein}. 

\begin{theorem}\label{th_red}{(Reduction theorem):}
Let $X$ be a complex manifold, $H$ a disc functional on $\A = \{ f\in \A_X ; f(0) \notin \sing(\omega)\}$ 
and $\omega$ a positive, 
closed $(1,1)$-current on $X$. 
The envelope $EH$ is $\omega$-plurisubharmonic if it satisfies the following.
\begin{enumerate}
 \item[(i)]\ $E\Phi^*H$ is $\Phi^*\omega$-plurisubharmonic for every holomorphic 
submersion $\Phi$ from a complex manifold where $\Phi^*\omega$ has a global potential
and for every $a \in \sing(\omega)$ we have $\limsup_{X\setminus \sing(\omega) \ni z\to a} EH(z) = EH(a)$.
 \item[(ii)]\ There is an open cover of $X$ by subsets $U$, with $\omega$-pluripolar 
 subsets $Z \subset U$ and local potentials $\psi$ on $U$, $\psi^{-1}(\{-\infty\}) \subset Z$, such that
 for every 
$h \in \A_U$, $h(\D) \nsubseteq Z$ the function $t \mapsto (H(h(t)) + \psi(h(t)))^\star$ is dominated by
an integrable function on $\T$.
 \item[(iii)]\ If $h \in \A_X$, $h(0) \notin \sing(\omega)$, 
 $t_0 \in \T \setminus h^{-1}(\sing(\omega))$ and $\epsilon > 0$,
then $t_0$ has a neighbourhood $U$ in $\C$ and there is a local potential $\psi$ in a neighbourhood of $h(U)$ 
such that for all sufficiently small arcs $J$ in $\T$ containing $t_0$ 
there is a holomorphic map
$F\colon D_r \times U \to X$ such that $F(0,\cdot) = h|_U$ and 
$$
\frac{1}{\sigma(J)} \int_J \big( H(F(\cdot,t))+\psi(F(0,t)) \big)^\star \ d\sigma(t) \leq (EH+\psi)(h(t_0)) + \epsilon.
$$
\end{enumerate}
\end{theorem}

Before proving the theorem we show that $H_{\omega,\phi}$ satisfies the conditions
\textsl{(i)-(iii)} and consequently that Theorem 1.1 follows from it.

\begin{prooftx}{Proof of Theorem 1.1}
	Condition \textsl{(i)} follows from Lemma \ref{submersion} since it implies
	$E\Phi^*H_{\omega,\phi} = EH_{\Phi^*\omega,\Phi^*\phi}$ and 
	by Theorem \ref{th_stein}, $EH_{\Phi^*\omega,\Phi^*\phi}$ is $\Phi^*\omega$-psh.

Condition \textsl{(ii)} follows from the fact that 
$H_{\omega,\phi}(h(t)) + \psi(h(t)) = \phi(h(t)) + \psi(h(t))$, which extends to an upper semicontinuous 
function on $\T$ and is thus dominated by a continuous function.

Assuming $h$ and $t_0$ as in condition \textsl{(iii)} and $\epsilon >0$, set $x=h(t_0)$ and let 
$f \in \A_X$ such that $f(0) = x$ and 
$H_{\omega,\phi}(f) \leq EH_{\omega,\phi}(x) + \epsilon/2$. 
By Lemma 2.3 in \cite{LarSig:1998} there is an open neighbourhood $V$ of $x$ in $X$, $r>1$
and a holomorphic function $\tilde F\colon D_r \times V \to X$ such that
$\tilde F(\cdot,x) = f$ on $D_r$ and $\tilde F(0,z) = z$ on $V$. Shrinking $V$ if necessary, we assume 
$\psi$ is a local potential for $\omega$ on $V$. Let $U = h^{-1}(V)$ and define 
$F\colon D_r \times U \to X$ by $F(s,t) =  \tilde F(s,h(t))$.
By the Riesz representation (\ref{riesz}),
\begin{equation}\label{riesz2}
\big(H_{\omega,\phi}(F(\cdot,t)) + \psi(F(0,t))\big)^\star = \int_\T (\phi + \psi)^\star \circ F(s,t)\, d\sigma(s).
\end{equation}
Since the integrand is usc on $D_r \times U$, then it is easily verified that (\ref{riesz2}) is an
usc function of $t$ on $U$. That allows us by shrinking $U$ to assume that
$$
\big(H_{\omega,\psi}(F(\cdot,t))+\psi(F(0,t))\big)^\star 
\leq H_{\omega,\phi}(F(\cdot,t_0)) + \psi(F(0,t_0)) + \frac \epsilon 2 
$$
for $t \in U$. Then by the definition of $f = F(\cdot,t_0)$
$$
\big(H_{\omega,\phi}(F(\cdot,t)) + \psi(F(0,t))\big)^\star < EH_{\omega,\phi}(x) + \psi(x) + \epsilon, \quad \text{on } U.
$$
Condition \textsl{(iii)} is then satisfied for all arcs $J$ in $\T \cap U$.

Finally, if $\mathcal F_{\phi,\omega} = \emptyset$ then the only function which
is both dominated by $\phi$ and satisfies the subaverage property is the
constant function $-\infty$. We know that $EH_{\phi,\omega} \leq \phi$ and
the proof of the reduction theorem gives the subaverage property, this
implies that $EH_{\phi,\omega} = -\infty$.
\end{prooftx}

\medskip

We now prove that $EH$ is $\omega$-psh if $H$ satisfies the three conditions in Theorem \ref{th_red}.
The main work is to show that $h^*EH$
satisfies the subaverage property of $h^*\omega$-sh functions for a given analytic
disc $h$. 
This implies by Proposition \ref{equi} that $EH$ is $\omega$-psh and that
concludes the proof of the reduction theorem.

\begin{lemma}
 Let $H$ be a disc functional on an $n$-dimensional complex manifold $X$ and $\omega$ a 
 positive, closed $(1,1)$-current on $X$. If the envelope $E\Phi^*H$ is 
$\Phi^*\omega$-usc for every holomorphic submersion $\Phi$ from an $(n+1)$-dimensional 
polydisc into $X$ and for every $a \in \sing(\omega)$ we have $\limsup_{X\setminus \sing(\omega) \ni z\to a} EH(z) = EH(a)$, then
$EH$ is $\omega$-usc.
\end{lemma}

\begin{proof}
 First,  
let $U$ be a coordinate polydisc in $X$ such that there is a potential $\psi$ for $\omega$ defined on $U$.
Define $\Phi\colon U \times \D \to U$ as the projection. By assumption and Lemma \ref{pullback}
$$
EH(z) + \psi(z) = EH(\Phi(z,t)) + \psi(\Phi(z,t)) \leq E\Phi^*H(z,t) + \psi(\Phi(z,t)) 
< +\infty, 
$$
for every $z \in U\setminus \sing(\omega)$ and $t\in \D$.
Then $EH < +\infty$ on $X\setminus \sing(\omega)$ and $EH + \psi$ is bounded above in some neighbourhood of 
every point in $\sing(\omega) \cap U$, where $\psi$ is any local potential for $\omega$.

Let $x \in X\setminus \sing(\omega)$ and $\beta > EH(x)+\psi(x)$. Assume $f \in \A_X$ is a 
holomorphic disc defined on $D_r$ such that $f(0)=x$ and 
$H(f) + \psi(x) < \beta$.
By using a theorem of Siu \cite{Siu:1976} it is shown
in the proof of Lemma 2.3 in \cite{LarSig:1998} that for $\tilde r \in ]0,r[$
there exists a neighbourhood $U$ of the
graph $\{t,f(t)\}$ in $D_r \times X$ and a biholomorphism 
$$
\Psi \colon U \to D_{\tilde r} \times \D^n
$$
such that $\Psi(t,f(t)) = (t,0)$.
Let $\pi\colon \C \times X \to X$ be the projection and define $\Phi = \pi \circ \Psi^{-1}$. 
Clearly $f = \Phi \circ \tilde f$, where $\tilde f \in A_{D_{\tilde r} \times \D^n}$ 
is the lifting $\tilde f(t) = (t,0)$. 
By assumption and the fact that 
$$
E\Phi^*H(x)+\psi(x) \leq \Phi^*H(\tilde f)+\psi(x) = H(f)+\psi(x) < \beta
$$ 
there is a neighbourhood $W$ of
$0 \in D_{\tilde r} \times \D^n$ such that 
$$ 
(E\Phi^*H +\Phi^*\psi)^\star < \beta, \quad \text{ on } W.
$$
Then for every $z$ in the open set $\Phi(W)$
$$
(EH(z)+\psi(z))^\star = (\Phi^*EH(\tilde z)+\psi(z))^\star \leq 
(E\Phi^*H(\tilde z)+\Phi^*\psi(\tilde z))^\star < \beta,
$$
where $\tilde z \in \Phi^{-1}(\{z\})$.
This along with the definition of $EH$ at $\sing(\omega)$ shows that the envelope is $\omega$-usc.
\end{proof}

Now we turn to the subaverage property of the envelope.
\begin{prooftx}{Proof of reduction theorem}
We have already shown that the envelope $EH$ is $\omega$-usc, so by Proposition \ref{equi} we only need to show that
for a local potential $\psi$ on an open set $U \subset X$ and every disc $h \in \A_U$ such that
$h(0) \notin \sing(\omega)$ we have
\begin{equation}\label{in0}
EH(h(0)) + \psi(h(0)) \leq \int_\T (EH \circ h + \psi\circ h)^\star \, d\sigma.
\end{equation}
Observe that this is automatically satisfied if $EH(h(0)) = -\infty$, so we may assume $EH(h(0))$ 
is finite.
It is sufficient to show that for every $\epsilon > 0$ and every continuous function $v\colon U \to \R$
with $v \geq (EH+\psi)^\star$, there exists $g \in \A_X$ such that $g(0) = h(0)$ and 
$$
H(g) + \psi(h(0)) \leq \int_\T v \circ h\, d\sigma + \epsilon.
$$
Then by definition of the envelope $EH(h(0)) + \psi(h(0)) \leq \int_\T v\circ h \, d\sigma + \epsilon$ 
for all $v$ and $\epsilon$, and (\ref{in0}) follows.

We assume that $h$ is holomorphic on $D_r$, $r>1$ and $h(0) \notin \sing(\omega)$. 
It is easily verified (see proof of Theorem 1.2 in \cite{LarSig:2003}) that a function satisfying 
the subaverage property for all holomorphic discs
in $X$ not lying in a pluripolar set $Z$ is plurisubharmonic not only on $X \setminus Z$ but on $X$.
We may therefore assume that $h(\overline \D) \nsubseteq Z$.

Note that $h(\T) \setminus \sing(\omega)$ is dense in $h(\T)$ by the subaverage property of 
$\psi\circ h$ and the fact that $h(0) \notin \sing(\omega)$. 
Therefore by a compactness argument along with property \textsl{(iii)} we can find a finite number of closed arcs 
$J_1,\ldots,J_m$ in $\T$, each contained in 
an open disc $U_j$ centered on $\T$ and holomorphic maps $F_j\colon D_s \times U_j \to X$, $s \in ]1,r[$ 
such that $F_j(0,\cdot) = h|_{U_j}$ and,
using the continuity of $v$, such that
\begin{equation} \label{in1}
\underline{\int_{J_j}} \Big( H(F_j(\cdot,t))+\psi (F(0,t)) \Big)^\star\, d\sigma(t) 
\leq \int_{J_j} v\circ h\, d\sigma + \frac{\epsilon}{4}\sigma(J_j).
\end{equation}
We may assume that the discs $U_j$ are relatively compact in $D_r$ and have mutually disjoint closure. 
By the continuity of $v$ and condition \textsl{(ii)} we may also assume that
\begin{equation}\label{in2}
\int_{\T \setminus \cup_j J_j} |v\circ h|\, d\sigma < \frac{\epsilon}{4}
\end{equation}
and
\begin{equation}\label{in3}
\overline {\int_{\T \setminus \cup_j J_j}} \Big(H(h(w))+\psi(h(w))\Big)^\star\, d\sigma(w) < \frac{\epsilon}{4}.
\end{equation}

We now embed the graph of $h$ in $\C^4 \times X$ as follows
$$
K_0 = \{ (w,0,0,0,h(w)) ; w \in \overline \D \}
$$
and the graphs of the $F_j$'s as
$$
K_j = \{ (w,z,0,0,F_j(z,w)) ; w \in J_j, z\in \overline \D \}.
$$
Let $\Phi\colon \C^4 \times X \to X$ denote the projection. This function restricted to a smaller subset 
will be our submersion.

What is needed to find the disc $g$ we are looking for is a Stein neighbourhood $V$ of
the compact set $K = \cup_{j=0}^m K_j$ in $\C^4 \times X$ where we can solve $d\eta =\Phi^*\omega$. 
Then we have by Lemma \ref{global_potential} a global potential for the pullback $\Phi^*\omega$ and
the $\Phi^*\omega$-plurisubharmonicity of $E\Phi^*H$, given by property \textsl{(i)}, then gives the existence of $g$.

For convenience we let $U_0 = D_r$ and $F_0(z,w) = h(z)$.
In \cite{LarSig:2003}, using Siu's theorem \cite{Siu:1976} and slightly shrinking the $U_j$'s and the $s$,
L\'arusson and Sigurdsson define for $j = 0,\ldots,m$ a 
biholomorphism $\Phi_j$ from $U_j \times D_s^{n+3}$ onto its image in $\C^4 \times X$ satisfying
$$
\Phi_j(w,z,0) = (w,z,0,0,F_j(z,w)), \quad w\in U_j, z \in D_s
$$
for $j=1,\ldots,m$, and
$$
\Phi_0(w,0) = (w,0,0,0,h(w)), \quad w \in D_r.
$$
The image of each $\Phi_j$ is therefore biholomorphic to a $4+n$ dimensional polydisc. 
These extensions of the graphs above are defined such that the first coordinate is the identity map.
This tells us that these images are mutually 
disjoint for $j\geq 1$ and that the intersection of the image of $\Phi_0$ and $\Phi_j$ is a subset of
$U_j \times \C^{3}\times X$.

As in \cite{LarSig:2003} we let $U_j''$ and $U_j'$ be discs cocentric with $U_j$ such that
$$
J_j \subset U_j'' \subset\subset U_j' \subset\subset U_j,
$$
and we assume our $\Phi_0$ is the $\Phi_0$ after the modification made in \cite{LarSig:2003} which are 
necessary to have all but the first coordinate of $\Phi_j^{-1} \circ \Phi_0$ close to the identity.
This modification which is done by precomposing $\Phi_0$ with a holomorphic map is necessary for constructing the 
Stein neighbourhood. Importantly for our purpose it does not change the first coordinate.

Now, for each $w \in \overline{U_j'}$ there is an $\epsilon_w >0$ such that 
$\Phi_0(w,D_{\epsilon_w}^{n+3}) \subset \Phi_j(w,D_s^{n+3})$. 
This holds by continuity for $\epsilon_w/2$ on a neighbourhood of $w$ in $U_j$.
By compactness of $\cup_{j=1}^m \overline{U'_j}$ there is an $\epsilon$ independent of 
$w$ such that $\Phi_0(w,D_{\epsilon}^{n+3}) \subset \Phi_j(w,D_s^{n+3})$ for $w \in \overline{U'_j}$.
We now restrict $\Phi_0$ to $D_r \times D_\epsilon^{n+3}$, then the intersection of the 
images $\Phi_0(D_r \times D_\epsilon^{n+3})$ and $\Phi_j(U'_j \times D_s^{n+3})$ is 
$\Phi_0(U'_j \times D_\epsilon^{n+3} )$

We define $V_j = \Phi_j(U_j \times D_s^{n+3})$, 
$V_0 = \Phi_0(U_0 \times D_\epsilon^{n+3})$ and $U = \cup_{j=0}^m V_j$.
To solve $d\eta = \Phi^*\omega$ on $U$ it then enough to show that the cohomology 
$H^2(U)$ is zero. This can be done using the exact 
Mayer-Vietoris sequence (\cite{Spivak:I} Ch.~11, Theorem 3), 
$$
\ldots \to H^q(M\cup N) \to H^q(M) \oplus H^q(N) \to H^q (M \cap N) \to H^{q+1}(M\cup N) \ldots.
$$
We start by letting $M=V_0$ and $N=V_1$, these sets and their intersection are biholomorphic to 
a polydisc, so they are smoothly contractable and then by Poincar\'e's lemma
$H^2(V_j)=H^1(V_0 \cap V_1)= 0$. 
Consequently, we see from the Mayer-Vietoris sequence
$$
\ldots \to H^1(V_0 \cap V_1) \to H^2(V_0 \cup V_1) \to H^2(V_0) \oplus H^2(V_1) \to \ldots
$$
that $H^2(V_0 \cup V_1) = 0$. Next we let $M=V_0 \cup V_1$ and $N = V_2$, then 
$H^1( (V_0\cup V_1) \cap V_2) = H^1( V_0 \cap V_2)=0$
since $V_1$ and $V_2$ are disjoint. The sequence above tells us then that $H^2(V_0 \cup V_1 \cup V_2) = 0$.
Iterating this process for all the $V_j$'s we finally see that $H^2(U) = H^2(\cup_{j=0}^m V_j) = 0$.

Next step is to find a Stein neighbourhood $V$ of $K$ which is a subset of $U$. 
Note that $K$ only relies on the holomorphic functions $h$ and $F_j$. Therefore $V$ can be constructed
in exactly the same way as in \cite{LarSig:2003}. It is done by defining a continuous strictly 
plurisubharmonic exhaustion function $\rho$ on $U$. This function is positive and satisfies
$K \subset \rho^{-1}[0,\frac 12]$. Finally $V$ is defined as $\rho^{-1}([0,1))$.

Then by Lemma \ref{global_potential} we have a global potential on $V$ for $\Phi^*\omega$.
By property \textsl{(i)} the envelope $E\Phi^*H$ is then $\Phi^*\omega$-psh on $V$
and if $\tilde h \colon D_r \to V$ is the lifting
$w\mapsto (w,0,0,0,h(w))$ of $h$, then 
$$
E\Phi^*H(\tilde h(0)) + \Phi^*\psi(\tilde h(0)) 
\leq \int_\T (E\Phi^*H \circ \tilde h + \Phi^*\psi \circ \tilde h)^\star \, d\sigma.
$$
Since $EH(h(0))\neq -\infty$ then $-\infty < \Phi^*EH(\tilde h(0))  \leq E\Phi^*H(\tilde h(0))$ and 
we may assume
there is a disc $\tilde g \in \A_V$ such that $\tilde g(0) = \tilde h(0)$ and 
$\Phi^*H(\tilde g) \leq E\Phi^*H(\tilde g(0)) + \epsilon/4$.
Define the disc $g = \Phi \circ \tilde g \in \A_X$, then $g(0) = h(0)$ and since $H(g)=\Phi^*H(\tilde g)$ 
and $\Phi^*\psi(\tilde h) = \psi(h)$, 
\begin{equation}\label{in4}
H(g) + \psi(h(0))  \leq 
\int_\T (E\Phi^*H \circ \tilde h + \psi\circ h)^\star\, d\sigma + \frac{\epsilon}{4}.
\end{equation}

For $w \in J_j$, $1\leq j\leq m$ we have a lifting of $F_j(\cdot,w)$ by $\Phi$ given by 
$z \mapsto (w,z,0,0,F_j(z,w))$.
Clearly $0 \mapsto \tilde h(w)$, so $E\Phi^*H(\tilde h(w)) \leq \Phi^*H(\tilde F_j(\cdot,w)) = H(F_j(\cdot,w))$. 
However, if $w \in \T \setminus \cup_j J_j$
then $E\Phi^*H(\tilde h(w)) \leq \Phi^*H(\tilde h(w)) = H(h(w))$. 
Therefore,
$$
\int_\T E\Phi^*H \circ \tilde h\, d\sigma 
\leq \sum_{j=1}^m \underline{\int_{J_j}} H(F_j(\cdot,w))\, d\sigma(w) 
+ \overline{\int_{\T\setminus \cup_j J_j}} H(h(w))\, d\sigma(w).
$$
Adding the integral of $\psi(h)$ to both sides of this inequality and using 
the inequalities (\ref{in1}) and (\ref{in3}) we see that
$$
\int_\T (E\Phi^*H \circ \tilde h + \psi\circ h)^\star\, d\sigma \leq 
\int_{\cup J_j} v\circ h\, d\sigma + \frac{\epsilon}{4}\sigma(\cup_j J_j) + \frac{\epsilon}{4}.
$$
Then by using first (\ref{in4}) and then (\ref{in2}) we have finally
$$
H(g) + \psi(h(0))  \leq \int_{\cup J_j} v\circ h\, d\sigma + \frac{3}{4}\epsilon
< \int_\T v\circ h\, d\sigma + \epsilon.
$$
\end{prooftx}

\bibliographystyle{siam}
\bibliography{bibref}

\bigskip
{\small
Science Institute, University of Iceland, Dunhaga 3, IS-107
Reykjavik, Iceland

E-mail: bsm@hi.is
}
\end{document}